\documentclass[11pt,reqno]{amsart}
\usepackage{amsmath,amsthm,amsfonts,amssymb,amscd,latexsym, calrsfs}
\usepackage[all]{xy}

\usepackage{amssymb}
\usepackage{mathrsfs}

\numberwithin{equation}{section}

\def\F{\mathbb{F}}

\def\su{{\subseteq}}
\def\la{{\langle}}
\def\ra{{\rangle}}

\def\J{\mathcal J}

\def\w{{\omega}}

\def\dim{{\rm dim\,}}
\def\char{{\rm char}}

\def\Proof{\noindent{\sl Proof.}\ }
\def\qed{{\hfill $\Box$ \medbreak}}

\newtheorem{defi}{Definition}[section]
\newtheorem{thm}[defi]{Theorem}
\newtheorem{lem}[defi]{Lemma}
\newtheorem{cor}[defi]{Corollary}

\newtheorem{eg}[defi]{Example}
\newtheorem{prop}[defi]{Proposition}
\newtheorem{rem}[defi]{Remark}

\usepackage[top=3.4cm,bottom=3.4cm,left=3.4cm,right=3.4cm, textwidth = 5in]{geometry}

\begin{document}

\title[Perfect and semiperfect restricted enveloping algebras]{Perfect and semiperfect restricted enveloping algebras}

\author{\textsc{Salvatore Siciliano}}
\address{Dipartimento di Matematica e Fisica ``Ennio De Giorgi", Universit\`{a} del Salento,
Via Provinciale Lecce--Arnesano, 73100--Lecce, Italy}
\email{salvatore.siciliano@unisalento.it}

\author{\textsc{Hamid Usefi}}
\address{Department of Mathematics and Statistics,
Memorial University of Newfoundland,
St. John's, NL,
Canada, 
A1C 5S7}
\email{usefi@mun.ca}

\thanks{The research of the second author was supported by NSERC of Canada under grant \# RGPIN 418201}

\begin{abstract}  For a restricted Lie algebra $L$,  the conditions under which its restricted enveloping algebra $u(L)$ is semiperfect are investigated. Moreover, it is proved that $u(L)$ is left (or right) perfect if and only if $L$ is finite-dimensional.
\end{abstract}

\date{\today}

  \maketitle

\section{Introduction}
Let $R$ be a ring with unity and denote by $\J(R)$ the Jacobson radical of $R$. We recall that $R$ is said to be \emph{semiperfect} if $R/\J(R)$ is Artinian and idempotents of $R/\J(R)$ can be lifted to $R$.  Semiperfect rings, introduced by H. Bass in \cite{Bass}, turn out to be a significant class of rings from the viewpoint of homological algebra and representation theory, since they are precisely the rings $R$ for which all finitely generated left or right $R$-modules have a projective cover (see e.g. \cite{L}, Chapter 8, \S 24). Clearly,  one sided Artinian rings and local rings  are semiperfect.

Recall that $R$ is called \emph{ left perfect} if all left $R$-modules have projective covers. Right perfect rings are defined in an analogous way.  
The pioneering work on perfect rings was carried out by H. Bass in 1960  and most of the main characterizations of these rings are contained in his celebrated paper \cite{Bass}. In particular, it follows from Bass' results  that the   following conditions  are equivalent: $R$ is left perfect;  every flat left $R$-module is projective; $R/\J (R)$ is Artinian and for every sequence $\{a_i\}$ in $\J (R)$ there exists an integer $n$ such that $a_1a_2\cdots a_n=0$; $R$ satisfies the descending chain condition on principal right ideals.
It should be mentioned that, while semiperfectness is a left-right symmetric property, there exist rings which are perfect on one side but not on the other (see \cite{Bass}, Example 5 on page 476). However, right and left perfectness are clearly equivalent conditions provided $R$ has a nontrivial involution. For instance, this is the case
when $R$ is a group algebra or an  (ordinary or restricted) enveloping algebra. 

Left perfect group rings were characterized by G. Renault in \cite{R} and, independently, by S. M. Woods in \cite{W1}. It turns out that, for a group $G$  and a field $\F$, the group algebra $\F G$ is left perfect if and only if $G$ is finite. Subsequently, a generalization of these results to semigroup rings has been carried out by J. Okni\'nski in \cite{O}. On the other hand, although semiperfect group algebras have been also investigated in several papers (see e.g. \cite{Bur, G, V, W2}),  a full characterization is not available yet. The best partial result in this direction  was obtained by J.M. Goursaud in \cite{G} where,  under the assumption  that the group $G$ is locally finite,  it is proved that $\F G$ is semiperfect if and only if either $\char \,\F=0$ and $G$ is finite or $\char \, \F=p>0$ and $G$ has a normal $p$-group of finite index.

In this paper, we consider these problems in the setting of enveloping algebras.  
For a restricted Lie algebra $L$ over a field $\F$ of characteristic $p>0$, we denote by $u(L)$ the restricted enveloping algebra of $L$ and by $L^\prime$ the derived subalgebra of $L$. Recall that a  subset $S$ of $L$ is said to be $p$-nil if every element of $S$ is $p$-nilpotent. We first provide some necessary and sufficient conditions on $L$ such that  $u(L)$ is semiperfect. We show that if $u(L)$ is semiperfect, then every element of $L$ is $p$-algebraic and $L$ does not contain any infinite-dimensional torus.  Under the assumptions that $\F$ is perfect 
 and $L^\prime$ is $p$-nil  we prove that for a locally finite-dimensional restricted Lie algebra $L$, $u(L)$ is semiperfect if and only if $L$ contains a $p$-nil restricted ideal of finite codimension.
 We construct an  infinite-dimensional  abelian restricted Lie algebra $L$ over an imperfect field $K$  such that $u(L)$ is semiperfect and $L$ has no nonzero $p$-nil restricted ideal. Hence, the perfectness assumption on the ground field is necessary.
It turns out that the  structure of semiperfect ordinary enveloping algebras $U(L)$ of arbitrary Lie algebras $L$ is  trivial and also quickly discussed.
    
In the last section, we prove that $u(L)$ is left (or right) perfect if and only if $L$ is finite-dimensional, which represents the  Lie-theoretic analogue of the aforementioned result of Renault and Woods.    

\section{Semiperfectness}

Throughout the paper, $\F$ denotes a field. Of course, every field  is  a left and right  perfect ring. To avoid any possible confusion, we recall that $\F$ is called a perfect field if every  finite extension of $\F$ is separable.

It would be an interesting problem to characterize perfect and semiperfect Hopf algebras. Our first result, which we actually use later in  the setting of restricted enveloping algebras,  is a step towards this goal. Recall that an element $\Lambda $ in a Hopf algebra $H$ with comultiplication $\Delta$ and  counit $\epsilon$ is called a left  integral if  $h \Lambda=\epsilon(h) \Lambda$, for every $h\in H$. Right integrals are defined analogously. Moreover, an element $x$ of $H$ is called \emph{group-like} if $x\neq 0$ and $\Delta(x)=x \otimes x$ and \emph{primitive} if $\Delta(x)=1\otimes x + x\otimes 1$.

\begin{thm}\label{hopf}
Suppose that a Hopf algebra $H$ is semiperfect. Then $H$ contains no infinite chain of semisimple Hopf subalgebras. Moreover, if $H$ is generated as an algebra  by its group-like and primitive elements, then all such elements are algebraic.   
\end{thm}
\Proof  For every subspace $V$ of $H$, we write $\bar{V}$ for the image of $V$ in $H/\J(H)$. By contradiction, suppose that $H$ contains an infinite chain of semisimple Hopf subalgebras 
$$
0\subsetneq H_1\subsetneq H_2 \subseteq \cdots \subsetneq H_n \subsetneq \cdots .
$$
 By Lemma 5.3.1 of \cite{DNR}, each Hopf algebra $H_i$ is finite-dimensional. Moreover,  by Theorem 2.2.1 and Corollary 2.2.4 of \cite{M}, left and right integrals of $H_i$ coincide and form a 1-dimensional ideal of $H_i$ which is not contained in $\ker \epsilon_i$, where  $\epsilon_i$ denotes the counit of $H_i$.  
Therefore we can find  $h_i\in \ker \epsilon_i$ such that $(1+h_i)x=0=x(1+h_i)$, for every $x\in \ker \epsilon_i$. Note that for every $i\leq j$ one has $(1+h_j)(1+h_i)=1+h_j=(1+h_i)(1+h_j)$.   As a consequence, we have a chain 
$$
\overline{H(1+h_1)}\supseteq \overline{H(1+h_2)} \supseteq \cdots \supseteq  \overline{H(1+h_n)} \supseteq \cdots
$$
of left ideals of $\bar{H}$. Now, since $H$ is semiperfect, the algebra $\bar{H}$ is left Artinian and so we must have 
$$
 \overline{H(1+h_n)}= \overline{H(1+h_{n+1})}
$$
for some $n$.  Therefore, as the images of $1+h_{n}$ and $1+h_{n+1}$ in $\bar{H}$ are commuting idempotents  generating the same left ideal, they must coincide. Thus we have $h_{n+1}-h_n \in \J(H)$. But we have
$$
(h_n-h_{n+1})^2=(h_{n}+1)^2+(h_{n+1}+1)^2-2(h_{n}+1)(h_{n+1}+1)=h_n-h_{n+1}.
$$
 We deduce that $h_{n+1}=h_n$, because  $\J(H)$ is free of nonzero idempotents. In particular, $1+h_{n+1}\in H_n$, a contradiction, yielding the first part of the theorem. 

Now, if  $H$ is generated as an algebra  by its group-like and primitive elements, then $H$ is clearly pointed. Let $x$ be a group-like or a primitive element of $H$ and suppose, by contradiction, that $x$ is trascendental.  Set $x_0=x$ and define recursively $x_{i+1}=x_i - x_i^2$ for every $n\geq 0$. Since $H/\J(H)$ is (right) Artinian, we deduce  from Lemma 3.1 in \cite{W2} that  there exists a positive integer $m$ such that $1-x_m$ has a right inverse $y$ in $H$. Let $B$ denote the Hopf subalgebra of $H$ generated by $x$, so that we have   $1-x_m\in B$.  Then $B$ is either the polynomial algebra $\F[x]$ (when $x$ is primitive) or the Laurent polynomial algebra  $\F[[x]]$  (when $x$ is group-like). In particular, $B$ has no nonzero zero-divisors. 
 Therefore, as $H$ is pointed, it follows from the main theorem in \cite{Rad}  that  $H$ is a free left $B$-module, say $H=\oplus_{h\in S} Bh$.  Note that by the proof of main  theorem in \cite{Rad}, we can assume that $1\in S$. Now we can deduce  
  that $y$ is contained in $B$. But the units of $\F[x]$ are non-zero elements of $\F$, and units of $\F[[x]]$ are of the form $\alpha x^{k}$, where $0\neq \alpha\in \F$ and $k$ is an integer. This yields  a contradiction, completing the proof.
\qed

 As group algebras are obviously generated by group-like elements,  we remark that Theorem \ref{hopf}  generalizes Theorem 3.2  in \cite{W2}.

Let $L$ be a restricted Lie algebra over a field of characteristic $p>0$. For $S\subseteq L$ we denote by $\langle S \rangle_p$ the restricted subalgebra generated by $S$. An element $x$ of $L$ is said to be \emph{$p$-algebraic} if $\dim \langle x \rangle_p<\infty$ and \emph{$p$-transcendental}, otherwise.   
As ordinary and restricted enveloping algebras are generated by primitive elements,  Theorem \ref{hopf} yields the following:  

\begin{cor}\label{torus}
Let $L$ be restricted Lie algebra over a field of characteristic $p>0$ such that $u(L)$ is semiperfect. Then $L$ is $p$-algebraic and contains no infinite-dimensional torus.
\end{cor}
\Proof In view of Theorem \ref{hopf}, $L$ is clearly $p$-algebraic. Now suppose, by contradiction, that $L$ contains an infinite dimensional torus $T$. Clearly, $T$ gives rise to an infinite chain 
 $$
0\subsetneq T_1\subsetneq T_2 \subseteq \cdots \subsetneq T_n \subsetneq \cdots 
$$
of finite-dimensional tori of $L$. Then, by \cite{H}, we see that $u(T_i)$ is a semisimple Hopf subalgebra of $u(L)$  for every $i$, which contradicts Theorem \ref{hopf}. \qed

Also, as the non-zero elements of  an arbitrary Lie algebra $L$ are transcendental in the ordinary enveloping algebra $U(L)$ of $L$,  from Theorem \ref{hopf} we deduce that $U(L)$ is semiperfect only in the trivial case:   

\begin{cor}\label{U(L)}
Let $L$ be a Lie algebra. Then the ordinary enveloping algebra $U(L)$ is semiperfect if and only if $L=0$.
\end{cor}

\begin{rem}\emph{ As the trivial module of an ordinary enveloping algebra $U(L)$ is not projective unless $L=0$, one can also prove Corollary \ref{U(L)} by using the fact that $U(L)$ is semiprimitive (see Theorem 8.3.14 in \cite{R1}) and so only projective $U(L)$-modules have projective covers (see e.g. \cite{R1}, Exercise 20 on page 326).}    
\end{rem}

Let $X$ be a basis of $L$ and $u\in u(L)$.
 Then, by the Poincar\'e-Birkhoff-Witt (PBW) Theorem for restricted Lie algebras (see \cite{SF}, Theorem 2.5.1), there exists a finite subset of $X$, denoted by $\text{Supp}(u)$, such that
 $u$ can be written as a linear combination of PBW monomials in the elements of $\text{Supp}(u)$ only. We recall that a restricted Lie algebra is said to be \emph{locally finite-dimensional} if all of its finitely generated restricted subalgebras are finite-dimensional.

\begin{prop}\label{nil}
Let $L$ be a locally finite-dimensional restricted Lie algebra over a field of characteristic $p>0$. Then 
$\mathcal{J}(u(L))$ is nil. In particular, $u(L)$ is semiperfect if and only if $u(L)/\mathcal{J}(u(L))$ is Artinian.
\end{prop}
\Proof Let $u\in \mathcal{J}(u(L))$. Since $\text{Supp}(u)$ is finite and $L$ is 
locally finite-dimensional, we have that $H=\la \text{Supp}(u)\ra_p$ is finite-dimensional. Since $u\in H$, we deduce that $u(L)$ is algebraic and the assertion follows from \cite[\S 1.10]{J}.
The second part is clear as idempotents can be always lifted modulo a nil ideal.\qed

We need the following result for later use, however we were unable to find a reference for it and as such we  include a proof   for completeness.

\begin{prop}\label{R/N}
Let $R$ be a ring and $N$ a nil ideal of $R$. If $R/N$ is semiperfect, then so is $R$.
\end{prop}
\Proof Note that $\mathcal{J} (R/N)=\mathcal{J} (R)/N$. Hence, we have
\begin{align*}
(R/N)/\mathcal{J} (R/N) \cong (R/N)/(\mathcal{J} (R)/N) \cong R/\mathcal{J} (R).
\end{align*}
Since $R/N$ is Artinian, it follows that $R/\mathcal{J} (R)$ is Artinian.

Now, let $e\in R$ such that the image of $e$ modulo $\mathcal{J} (R)$ is an idempotent in $R/\mathcal{J} (R)$. Since $R/N$ is semiperfect, there exists $f\in R$ such that the image of $f$ modulo $N$ is an idempotent in $R/N$ and $e-f\in \mathcal{J} (R)$. But $N$ is nil, and it is well known that $f$ can be lifted to an idempotent of $R$, that is, there exists $g\in R$ such that $f-g\in N\su \mathcal{J} (R)$. Notice that $e-g\in  \mathcal{J} (R)$ and we are done. \qed

\begin{thm}\label{locally finite}
Let $L$ be a locally finite-dimensional restricted Lie algebra over a field of characteristic $p>0$.
If $L$ has a   $p$-nil restricted ideal of finite codimension, then $u(L)$ is semiperfect.
\end{thm}
\Proof Let $P$ be a $p$-nil restricted ideal of finite codimension. We first prove that the ideal $I=Pu(L)$ is nil. To do so, let $u\in I$. Then $u=\sum_i^{n}p_iu_i$, where each $p_i\in P$ and each $u_i\in u(L)$. Let $\mathcal{B}$ be a basis of $L$ containing a basis of $P$. Note that there exists a finite subset $S$ of $\mathcal{B}$ such that all the elements $p_i, u_i$ are in $u(H)$, where $H$ is the restricted subalgebra of $L$ generated by $S$. Now, let $Q$ be the restricted ideal of $H$ generated by all the $p_i$'s. Then $Q\subseteq P\cap H$, and so we see that $Q$ is a finite-dimensional $p$-nilpotent restricted ideal of $H$. Note that $Qu(H)$ is associative nilpotent.  Since $u\in Qu(H)$, 
 $u$ is nilpotent and as such  $I$ is nil. 
Now, note that 
\begin{align*}
u(L/P)\cong u(L)/I.
\end{align*}
Since $L/P$ is finite-dimensional, clearly $u(L/P)$ is semiperfect. 
Hence, $u(L)/I$ is semiperfect and it follows from Proposition \ref{R/N} that $u(L)$ is semiperfect.\qed

 We remark that an argument similar to the one used in Theorem \ref{locally finite} yields an alternative  and shorter proof of the analogous result for group rings obtained  in  \cite{W2}.

For a group algebra $\F G$ of a locally finite group $G$ over a field $\F$ of characteristic $p>0$, it is shown by J.M. Goursaud  (see \cite{G}, Th\'{e}or\`{e}me 8)  that $\F G$ is semiperfect if and only $[G:\mathcal{O}_p(G)]<\infty$, where $\mathcal{O}_p(G)$ is the largest normal $p$-subgroup of $G$.  In analogy, we let $\mathcal{O}_p(L)$ be the sum of all $p$-nil restricted  ideals of $L$. Since the sum of every two $p$-nil restricted ideals of $L$ is again $p$-nil, it follows that $\mathcal{O}_p(L)$ is also a $p$-nil restricted ideal. Clearly, if $L$ is finite-dimensional then $\mathcal{O}_p(L)$ is just $\mathrm{rad}_p(L)$, the $p$-radical  of $L$.   
Therefore one might expect that the analogous result holds for restricted Lie algebras, that is,  
$u(L)$ is semiperfect if and only $\dim L/\mathcal{O}_p(L)<\infty$. In Proposition \ref{semiperfect} we provide a partial converse of Theorem \ref{locally finite} and, on the other hand,   in Example \ref{perfect field}  we 
present an infinite-dimensional  abelian restricted Lie algebra $L$ over an imperfect field   such that $u(L)$ is semiperfect and  yet $\mathcal{O}_p(L)=0$.  Thus, even in the abelian case, the Lie-theoretic analogue of Goursaud's Theorem does not hold in general.

\begin{prop}\label{semiperfect}
 Let $L$ be a locally finite-dimensional restricted Lie algebra over a perfect field $\F$ of characteristic $p>0$. If $u(L)$ is semiperfect and $L'$ is $p$-nil,  then $\mathcal{O}_p(L)$ has finite codimension in $L$.
\end{prop}
\Proof First we show that the associative ideal  $I=\mathcal{O}_p(L)u(L)$ is nil.
Let $u\in I$. Then $u=\sum x_i u_i$,  where $x_i\in \mathcal{O}_p(L)$ and $u_i\in u(L)$. Now let $H$ be  the restricted subalgebra generated by all the $x_i$'s and all the $\text{Supp}(u_i)$'s. Note that $H$ is finite-dimensional.
Let $N$ be the restricted ideal of $H$ generated by the $x_i$'s. 
Since $N\su \mathcal{O}_p(L)$, we deduce that $N$ is finite-dimensional and $p$-nilpotent. It follows that $Nu(H)$ is nilpotent. In particular, as $u\in Nu(H)$, we deduce that $u$ is nilpotent. Therefore, $I$ is a nil ideal and as such $I\su \J (u(L))$. We now claim that $I= \J (u(L))$.
Let $\mathcal{L}=L/\mathcal{O}_p(L)$. Since $L^\prime\su \mathcal{O}_p(L)$,  $\mathcal{L}$ is  free of nonzero $p$-nilpotent elements.  We claim that $u(\mathcal{L})\cong u(L)/I$ is reduced. Let $u$ be a nilpotent element of $u(\mathcal{L})$. There  exist  elements $x_1,\ldots,x_n$ in $L$ that are linearly independent modulo $\mathcal{O}_p(L)$ such that
$$
u=\sum \alpha x_1^{i_1}\cdots x_n^{i_n} \quad \text{ modulo } I.
$$
Thus, $u^{p^r}=0$  for some positive integer $r$. Hence, by the PBW Theorem,  the elements $x_1^{[p]^r},\ldots,x_n^{[p]^r}$ are linearly dependent modulo $\mathcal{O}_p(L)$. Let $\beta_1,\beta_2,\ldots,\beta_n\in \F$, not all zero, such that
$$
\sum \beta_k x_k^{[p]^r}\in \mathcal{O}_p(L).
$$
Since $\F$ is a perfect field, we get
$$
\left(\sum \beta_k^{\frac{1}{p^r}} x_k\right)^{[p]^r} \in \mathcal{O}_p(L).
$$
But $\mathcal{L}$ has no nonzero $p$-nilpotent elements. We deduce that
$\sum \beta_k^{\frac{1}{p^r}} x_k \in \mathcal{O}_p(L)$, which contradicts the fact that 
 $x_1,\ldots,x_n$ are linearly independent modulo $\mathcal{O}_p(L)$. Hence, $u(\mathcal{L})$ is reduced.  On the other hand,  by Proposition \ref{nil}, $\J(u(\mathcal{L}))$ is nil.  We deduce  that $\J(u(\mathcal{L}))=0$. Thus, $\J (u(L))\subseteq I$, as claimed. At this stage, as $u(\mathcal{L}) \cong  u(L)/I =u(L)/\J (u(L))$ is Artinian, we deduce from  \cite{LZ}  that $\mathcal{L}$ is indeed finite-dimensional. Therefore $\mathcal{O}_p(L)$ has finite codimension in $L$, as required. \qed 
 
As an immediate consequence of Theorem \ref{locally finite}  and Proposition \ref{semiperfect} we have the following characterization of semiperfect commutative restricted enveloping algebras over perfect fields:
 
\begin{cor}\label{s.p.}
Let $L$ be an abelian restricted Lie algebra over a perfect field of characteristic $p>0$. Then $u(L)$ is semiperfect if and only if $\mathcal{O}_p(L)$ has finite codimension in $L$.
\end{cor}

It is worth mentioning that, in general, the properties of perfectness or semiperfectness of an algebra are not preserved under extensions of the ground field. For instance, consider the following:

\begin{eg}\label{field} \emph{
Let $\F$ be an algebraically closed field and let $E=\F(t,t^{\frac{1}{2}}, t^{\frac{1}{4}},\ldots )$ be the extension of the field $\F(t)$ of rational functions in the variable $t$ obtained by adjoining all the 2-power roots of $t$. Consider the automorphism $\alpha: E \longrightarrow E$, $f(t)\mapsto f(t^2)$ and form the division algebra $D=E((x;\alpha))$ of skew Laurent series  (for details, see e.g., Section 2.3 of \cite{C}). Obviously, $D$ is left perfect. Now, let $K=\F((x;\alpha))$. Then $K$ is a maximal subfield of $D$ and, clearly, the algebra $D$ is not algebraic over $K$. Moreover, as $\F$ is algebraically closed,  $K$ is a transcendental field extension of $\F$.  Therefore,  Theorem 2.4 of \cite{JP} assures that $\J(D\otimes_{\F}K)=0$. Finally, by Exercise 6 in Section IX.6 of \cite{Hu}, we see that  $D\otimes_{\F}K$ is not left Artinian and, indeed,  not semiperfect.    
}
\end{eg}
 
We remark that the analogue of the aforementioned result of Goursaud  for restricted Lie algebras fails without the 
perfectness assumption on the field.
We show as follows that the restriction on the ground field is  required. 
For simplicity, we construct our example over a field  $\F$  of characteristic 2
but this example could be extended to any positive characteristic.
In the following we denote by $K$ the field of rational functions in infinitely many indeterminates $t_1, t_2, \dots$ over $\F$, that is,    $K=\F(t_1 , t_2, \cdots)$.

\begin{thm}\label{perfect field}
 There exists an infinite-dimensional  abelian restricted Lie algebra $L$ over $K$  such that $u(L)$ is semiperfect and  $\mathcal{O}_p(L)=0$. 
 \end{thm}
\Proof Let $L$ be the abelian restricted Lie algebra over $K$ with basis $x, y_1, y_2, \ldots$ and power mapping defined by $x^{[2]}=x, y_i^{[2]}=t_i x$ for $i\geq 1$.  

We have that $\mathcal{O}_p(L)=0$. Indeed, let $z\in L$ such that $z^{[2]}=0$. Then $z=\alpha x+\sum \beta_i y_i$ for some $\alpha, \beta_i\in K$. We have
$$
z^{[2]}=(\alpha^2+\sum \beta_i^2 t_i )x=0.
$$
Hence, $\alpha^2+\sum \beta_i^2 t_i =0$. It is  easy to deduce that  $\alpha=\beta_i=0$ for all $i$. Hence, $z=0$ and $\mathcal{O}_p(L)=0$. We denote by $\bar u$ the image of an element $u$ modulo $\J(u(L))$.  
Note that $L$ is locally finite-dimensional. Hence, to show that $u(L)$ is semiperfect,  by the second part of Proposition \ref{nil} it is enough to prove that
$\bar R=u(L)/\J(u(L))$ is Artinian. For this, we will first show that $u(L)/\J(u(L))$ satisfies the descending chain condition on principal ideals by proceeding in the following steps.

\textbf{Claim 1:} \emph{Let $u\in u(L)$. Then $\bar u$ is invertible in $u(L)/\J(u(L))$ if and only if $u$ is  invertible in $u(L)$.}
Suppose $\bar u \bar v= \bar 1$, for some $v\in u(L)$.  Then $1-uv$ is nilpotent. Thus, $uv$ is invertible which implies that $u$ is invertible. 

\textbf{Claim 2:} \emph{Let $u=a+bx\in u(L)$, where $a, b\in K$. Then $u$ is invertible if and only if $a\neq 0$ and $a\neq b$.}
Suppose that  $a\neq 0$ and $a\neq b$. Then we have 
$$
(a+bx)(a^{-1}+b(a^2+ab)^{-1}x)=1.
$$
The converse follows from the fact that $x(x+1)=0$.

\textbf{Claim 3:} \emph{Let $0\neq d\in K^{2^{r}}( t_1^{2^{(r-1)}}, t_2^{2^{(r-1)}}, \ldots )$, where $r\geq 1$. Then there exists $w\in u(L)$ such that  $w^{2^r}=d^{-1} x$.} Suppose that 
$$
d= \beta^{2^r}  +\sum_{n\geq 1} \gamma_{i_1, \ldots, i_n}^{2^r} (t_{i_1}\cdots t_{i_n})^{2^{(r-1)}}
$$
for some $\beta$ and $\gamma_{i_1, \ldots, i_n}\in K$.
Put

\[
w=\frac{\beta x  +\sum_{n\geq 1} \gamma_{i_1, \ldots, i_n} y_{i_1}\cdots y_{i_n}}{\beta^2  +\sum_{n\geq 1} \gamma_{i_1, \ldots, i_n}^2 t_{i_1}\cdots t_{i_n}}.
\]
We observe that  $w^{2^r}=d^{-1} x$.

\textbf{Claim 4:} \emph{For any chain of principal ideals $0\subsetneq (\bar u_2)\su (\bar u_1)\subsetneq \bar R$, we have that $ (\bar u_2)= (\bar u_1)$.}
 Let $v\in u(L)$ such that   $\bar u_2=\bar u_1 \bar v$. We may assume that $\bar v$ is not invertible.
 We show that there exists $w\in u(L)$ such that $u_1-u_1vw$ is nilpotent.
Note that  $(y_i(x-1))^2=0$,  for every $y_i$. It follows that $\w(L)(x-1)\su \J(u(L))$.
In particular, by the PBW Theorem, every element $u\in u(L)$ can be written in the form 
\begin{align*}
u=\alpha + \beta x +\sum \gamma_{i_1, \ldots, i_n} y_{i_1}\cdots y_{i_n} \quad  \text{ modulo } \J(u(L)). 
\end{align*}
Thus, for  large  $r$, we have 
$$
u^{2^r}=\alpha^{2^r} + \beta^{2^r} x +\sum \gamma_{i_1, \ldots, i_n}^{2^r} (t_{i_1}\cdots t_{i_n})^{2^{(r-1)}}x.
$$
We deduce that  $u^{2^r}=a + bx$, where $a\in K^{2^r}$ and $b\in K^{2^{r}}( t_1^{2^{(r-1)}}, t_2^{2^{(r-1)}}, \ldots )$.

Let $u_1^{2^r}=a + bx$ and $v^{2^r}=c + dx$, where $a, c\in K^{2^r}$ and $b, d\in K^{2^{r}}( t_1^{2^{(r-1)}}, t_2^{2^{(r-1)}}, \ldots )$. Since $\bar u_1$ is not invertible, it follows from Claim 2 that either  $a= 0$ or $a= b$. Similarly, since $\bar v$ is not invertible,  either  $c= 0$ or $c=d$.

Suppose first  that $a=0$. Then, if $c=d$, we get $\bar u_2^{2^r}=\bar u_1^{2^r} \bar v^{2^r}=bcx(1+x)=0$.
Hence, $u_2\in \J(u(L))$ and so $\bar u_2=0$ which is a contradiction to the assumptions. We deduce that, 
if $a=0$, then $c=0$ and $c\neq d$.  Now, by Claim 3, there exists $w\in u(L)$ such that  $w^{2^r}=d^{-1}x$.
It is easy to see that $u_1^{2^r}-(u_1vw)^{2^r}=0$.
Now, assume that $a\neq 0$ and so $a=b$. In a similar way as above, we can see that  $c\neq 0$ and $c=d$.
  Since $c\in K^{2^r}$, there exists $w\in u(L)$ such that $w^{2^r}=c^{-1}$.
  It is clear that $u_1^{2^r}-(u_1vw)^{2^r}=0$. This completes the proof of Claim 4.

It now follows from Claim 4 that  $u(L)/\J(u(L))$ satisfies the descending chain condition on principal ideals. 
Since $u(L)/\J(u(L))$ is semiprime, we deduce from Theorem 10.24 in Chapter 4 of \cite{L}  that 
 $u(L)/\J(u(L))$ is Artinian, as required.
 \qed

\section{Perfectness}
In this section, we characterize left perfect restricted enveloping algebras $u(L)$. As the antipode of any cocommutative Hopf algebra is an involution,  left and right perfectness are, in fact, equivalent conditions for $u(L)$. 
 
\begin{lem}\label{fg}
Let $L$ be a restricted Lie algebra over a field of characteristic $p>0$. If $u(L)$ is left perfect, then 
$L$ is finitely generated.
\end{lem}
\Proof Let $R=u(L)$. 
Suppose, by contradiction, that $L$ is not finitely generated and let $x_1, x_2,\ldots \in L$ so that $x_{k+1} \notin \la  x_1, \ldots,  x_k \ra_p$.
Now consider the  principal right ideals
\begin{align*}
x_1R\supseteq x_1x_2R\supseteq \cdots \supseteq x_1x_2\cdots x_kR\supseteq\cdots .
\end{align*}
 Since, according to Bass' Theorem (see \cite{Bass}, Theorem P), $R$ satisfies the descending chain condition
on principal right ideals, there exists an integer $i$ such that  
\begin{align*}
x_1\cdots x_iR=x_1\cdots x_ix_{i+1}R.
\end{align*}
Therefore one has 
\begin{align}\label{x1xi}
x_1\cdots x_i-x_1\cdots x_ix_{i+1}v=0,
\end{align}
 for some $v\in u(L)$.  Let $H$ be the restricted subalgebra of 
 $L$ generated by $x_1, \ldots, x_i$.  We extend $x_1, \ldots, x_i$ to a basis $X$ of $H$ and extend $X$ 
 to a basis $X\cup Y$ of $L$. Since $x_{i+1}\notin H$, we can assume that $x_{i+1}\in Y$. It is well known that $u(L)$ is a free left $u(H)$-module, that is, 
 $$
 u(L)=\oplus_i u(H) u_i,
 $$
  where the $u_i$'s are distinct PBW monomials in the elements of $Y$.
Now, we observe that  Equation \eqref{x1xi} is not possible because  some monomials in the  PBW representation of $x_{i+1}v$ involve $x_{i+1}$ and this contradicts $u(L)$ being a free $u(H)$-module. 
Hence, $L$ must be finitely generated.
\qed

\begin{lem}\label{subalgebra}
Let $L$ be a restricted Lie algebra over a field of characteristic $p>0$. If $u(L)$ is left perfect, then so is 
$u(H)$  for every restricted subalgebra $H$ of $L$.
\end{lem}
\Proof 
Let $Y$ be a basis for a complementary vector subspace of $H$ in $L$. Since $u(L)$ is a free left $u(H)$-module, we have 
 \begin{equation}\label{subalg}
 u(L)=\oplus_k u(H) u_k,
 \end{equation}
  where the elements $u_k$ are distinct PBW monomials in the elements of $Y$. Let $h\in u(H)$ and let $I=hu(H)$ be the right ideal of $u(H)$ generated by $h$. Then (\ref{subalg}) implies that 
  $$
  J=\oplus_k I u_k
  $$  
is a principal right ideal of $u(L)$. As a consequence, every descending chain $\{I_k\}_k$ of principal right ideals of $u(H)$ gives rise to a descending chain  of principal right ideals of $u(L)$ which must stabilize, as $u(L)$ is left perfect. Therefore, by (\ref{subalg}), the chain $\{I_k\}_k$ stabilizes, which finishes the proof. \qed

 For a restricted Lie algebra $L$, we define
$$
\Delta(L)=\{x\in L \vert\, \dim [L,x]< \infty\}.
$$
Then $\Delta(L)$ is clearly a restricted ideal of $L$ which is the Lie algebra analogue of the FC-center of a group.

\begin{prop}\label{torsionFC}
Let $L$ be a restricted Lie algebra over a field of characteristic $p>0$. If $L=\Delta(L)$ and every element of $L$ is $p$-algebraic, then $L$ is locally finite-dimensional.
\end{prop}
\Proof Let $x_1,\ldots,x_n\in L$ and let $H$ be the (ordinary) subalgebra of $L$ generated by these elements. By Lemma 1.3 of \cite{BP} we have that $\dim [H,H]<\infty $. Since $H$ is spanned as a vector subspace by $x_1,\ldots,x_n$ and their Lie commutators, we see that $H$ is indeed finite-dimensional. Since  every element of $H$ is $p$-algebraic, it follows from Proposition 1.3(1) in Chapter 2 of \cite{SF} that $\langle x_1,\ldots,x_n \rangle_p$  is finite-dimensional, as required. \qed

The following result is an immediate consequence of Corollary 6.4  in \cite{BP}:

\begin{lem}\label{u(Delta)} Let $L$ be a restricted Lie algebra over a field of characteristic $p>0$. Then $u(L)$ is semiprime if and only if $u(\Delta(L))$ is semiprime.
\end{lem}

We can now prove the main result of this section:

\begin{thm}\label{perfect}
Let $L$ be a restricted Lie algebra over a field of characteristic $p>0$. Then $u(L)$ is left perfect if and only if $L$ is finite-dimensional.
\end{thm}
\Proof
One implication is clear. Suppose now that $u(L)$ is left perfect. By Lemma \ref{subalgebra}, $u(\Delta(L))$ is also left perfect and so, by Lemma \ref{fg}, $\Delta(L)$ is finitely generated. Note that, by Corollary \ref{torus},   every element of $\Delta(L)$ is $p$-algebraic. Hence,    by Proposition \ref{torsionFC}, we conclude that $\Delta(L)$ is finite-dimensional. Thus, it is enough to prove that $L/\Delta(L)$ is 
finite-dimensional.
Note that,  by Corollary 24.19 of \cite{L}, the ring $u(L/\Delta(L))\cong u(L)/\Delta(L)u(L)$ is left perfect, as well. Hence, we can replace $L$ by $L/\Delta(L)$. By   \cite[Lemma 2.7(i)]{RS}, we have that  $\Delta(L/\Delta(L))=0$. Therefore we can  assume that $\Delta(L)=0$.  Now,   Lemma \ref{u(Delta)} entails that $u(L)$ is semiprime and, by Bass' Theorem (see \cite{Bass}, Theorem P), $u(L)$ satisfies the descending chain condition on principal right  ideals.  But then, by Theorem 10.24 in \cite{L}, $u(L)$ is right Artinian. Since, by \cite{LZ},  right Artinian Hopf algebras are finite-dimensional, we conclude that $L$ is finite-dimensional. This completes the proof. \qed    

Some immediate consequences of Theorem \ref{perfect} are now in order. 
As the category of (left) $u(L)$-modules is equivalent to the category of restricted (left) $L$-modules, a combination of Theorem \ref{perfect} and  Theorem 24.25 in \cite{L} yields  the following:

\begin{cor}
Let $L$ be a restricted Lie algebra. The following conditions are equivalent:
\begin{enumerate}
\item every restricted flat $L$-module is projective;
\item $L$ is finite-dimensional.
\end{enumerate}
\end{cor}

Let $R$ be a ring and $\theta: G\to \text{Aut}( R)$  be a group homomorphism.  Then J.K. Park in \cite{park} proved that the smash product  $R \# \F G$ is left perfect if and only if $R$ is left perfect and $G$ is finite. In particular, in view of Theorem \ref{perfect}, if $R=u(L)$ for a restricted Lie algebra $L$, we deduce the following:

\begin{cor}
Let $A=u(L)\# \F G$. Then $A$ is left perfect if and only if $A$ is finite-dimensional.
\end{cor}

Note also that Theorem \ref{perfect} generalizes Corollary 6.6 in \cite{BP}.
Finally,  the celebrated structure theorem of Cartier-Kostant-Milnor-Moore (see e.g. \cite[\S 5.6]{M}) implies that every cocommutative Hopf algebra over an algebraically closed field of characteristic
zero can be presented as a smash product of a group algebra and an enveloping algebra.
So,  Corollary \ref{U(L)}  along with  Park's result  in \cite{park}  yields: 

\begin{cor}
Let $H$ be a cocommutative Hopf algebra over an algebraically closed field $\F$ of characteristic zero. Then $H$ is left perfect if and only if $H$ is the group algebra of a finite group $G$ over $\F$.
\end{cor}

\end{document}